\documentclass[10 pt]{amsart}
\usepackage{amscd,amsmath,amsthm,amssymb}
\usepackage{enumerate,varioref}
\usepackage{epsfig}
\usepackage{graphicx}
\newtheorem{thm}{Theorem}

\newtheorem{lem}[thm]{Lemma}
\newtheorem{prop}[thm]{Proposition}
\theoremstyle{definition}
\newtheorem{defn}[thm]{Definition}
\theoremstyle{remark}
\newtheorem{ex}[thm]{Example}
\newtheorem{rem}[thm]{Remark}

\numberwithin{equation}{subsection}
\newcommand{\R}{\mathbb{R}}

\newcommand{\C}{\mathbb{C}}

\newcommand{\ad}{a^{\dagger }}
\newcommand{\adsq}{a^{\dagger 2}}
\newcommand{\adcu}{a^{\dagger 3}}
\newcommand{\adqu}{a^{\dagger 4}}
\newcommand{\adn}{a^{\dagger n}}
\begin{document}

\title{Explicit Calculations for Anharmonic Oscillators Using Lie Algebras}
\author{Clark Alexander}

\email{gcalex@temple.edu} \maketitle \tableofcontents

\section{Algebraic Preliminaries}
The goal of the present paper is primarily to exhibit the effectiveness of using Lie algebras to 
compute explicit perturbation eigenvalues for quantum anharmonic oscillators in one dimension.  
There are however, several goals secondary in stature, but which merit discussion.  The first of these 
is to enable the reader to work with Weyl algebras in the abstract.  Two presentations of it arise 
readily in quantum physics.  In particular the first is the algebra of position and momentum 
operators in nonrelativistic mechanics wherein $[x,p]=i\hbar$.  The second is the algebra of 
ladder operators $[a,a^{\dagger}]=1.$  It is the second presentation with which we will be 
primarily concerned in this paper.  Of course, the first presentation may be made to look like 
the second by considering not $p$, but instead the simple derivative $\frac{d}{dx}$ whereby one has $[\frac{d}{dx},x]=1$.

\subsection{Normal Ordering and Weyl Binomial Coefficients}
For any abstract Weyl algebra determined by two elements $A$ and $B$ obeying $[A,B]=1$, an ordering of a polynomial in 
$A$ and $B$ will be said to be \emph{normally ordered} if all powers of $B$ appear to the left of powers
of $A$.  For example $A^3B^2$ is not normally ordered, but $B^2A^3$ is.  In our case $\ad$ will always be placed to the left of $a$.  
As it is well known in elementary quantum mechanics one may move back and forth between presentations of problems
in position-momentum coordinates and annihilator-creator coordinates with the following equivalences
\begin{eqnarray}
x = \frac{a+\ad}{\sqrt{2}}\\
p = \frac{a-\ad}{i\sqrt{2}}.\nonumber
\end{eqnarray}
Since this paper is concerned with anharmonic oscillators we will be concerned with $x^n$ in the potential.  
Thus, we need an efficient way of normally ordering $(a+\ad)^n$.
\begin{lem}
Let $A$ and $B$ determine a Weyl algebra so that $[A,B]=1$.  The normal ordering of $(A+B)^n$ is given by
\begin{equation}
(A+B)^n = \sum_{m=0}^n \sum_{k=0}^{\min\{n,n-m\}}\left\{\begin{array}{c}n\\m\end{array}\right\}_k B^{m-k}A^{n-m-k},
\end{equation}
where 
\begin{equation}
\left\{\begin{array}{c}n\\m\end{array}\right\}_k = \frac{n!}{2^k k! (m-k)!(n-m-k)!}
\end{equation}
is the \emph{Weyl binomial coefficient}.
\end{lem}

The proof of this lemma involves nothing more than counting commutations.
\begin{ex}
We will use the fourth order relation explicitly later, so here is an example of how the
Weyl coefficients factor in.
\[ 
(a+\ad)^4 =\adqu+4\adcu a+6\adsq a^2+4\ad a^3 + a^4 + 6\adsq +12\ad a + 6a^2 + 3.
\]
\end{ex}

\begin{rem}
Notice that 
\[
\left\{\begin{array}{c}n\\m\end{array}\right\}_k =
\left\{\begin{array}{c}n\\n-m\end{array}\right\}_k.
\]
\end{rem}

If one wishes to attempt calculations within a Weyl algebra it may be useful to compute with abstract elements $A,B$ first and then 
plug into a specific situation one has in mind.  One other useful tip is that if one has 
Weyl variables $A,B$ then it can be convenient to consider representing the algebra as $\frac{d}{dB},B$ or $A,-\frac{d}{dA}$.  
This becomes consistent with the first presentation considered.  For example, one peculiar formula which is arguably 
easier to compute with abstract Weyl variables is
\begin{equation}
\mu^{x\partial_x}f(x) = f(\mu x).
\end{equation}

\subsection{Baker-Campbell-Hausdorff and the Hadamard Lemma}
We will be concerned throughout much of this paper with exponentiating noncommuting variables.  
We run into a stopping block in trying to compute the exponentials explicitly.
The main issue is that for noncommuting variables $X,Y$ we see
\[
e^Ye^X \neq e^{Y+X} = e^{X+Y} \neq e^Xe^Y.
\]
The \emph{Baker-Campbell-Hausdorff} formula is the solution to $Z=\log(e^Xe^Y)$.  The explicit solution is formally given 
as symmetric sums and differences of nested commutators in $X$ and $Y$.  One may find this expression in 
nearly any textbook on advanced quantum mechanics.  We shall not be concerned, however, with isolated exponentials,
but rather expressions of the form
\[
e^X Y e^{-X}.
\]
Using elementary combinatorics and the Baker-Campbell-Hausdorff formula one can arrive at the \emph{Hadamard lemma}.
\begin{lem}
Let $X,Y$ be noncommuting variables then one has
\begin{equation}
e^X Y e^{-X} = Y + [X,Y] + \frac{1}{2!}[X,[X,Y]] + \frac{1}{3!}[X,[X,[X,Y]]] + \cdots
\end{equation}
If we allow the notation $[X^{(n)},Y] = [X,[\dots,[X,Y]]]$ then we may write more succinctly
\begin{equation}
e^X Y e^{-X} = \sum_{k=0}^{\infty}\frac{1}{k!}[X^{(k)},Y].
\end{equation}
\end{lem}

\subsection{The Formula Often Desired and Rarely Known}
One final assertion about Weyl variables in the algebraic preliminaries must be the formula 
\begin{equation}
[A^n,B^m] = \sum_{k=1}^{\min\{n,m\}}k!\left(\begin{array}{c}m\\k\end{array}\right)
\left(\begin{array}{c}n\\k\end{array}\right)B^{m-k}A^{n-k}.
\end{equation}

This formula is often left as an exercise in quantum mechanics texts and sometimes in homological algebra, 
but rarely is it completed.  One might jokingly say it is similar to the snake lemma in that no one knows
if it's really true since the only persons who have ever proven it are graduate students.  
All kidding aside, this is indeed the correct formula for normally ordering variables obeying the Weyl relation.

\section{The Method of [JA]}

The main impetus for this research comes from the paper [JA].  The goal of this section is to explicate in 
a reasonably clear manner the content of that paper.

The premise upon which [JA] begins is the idea that we can create a new Lie algebra from simply taking commutators 
of the unperturbed Hamiltonian $H_0$ and the new anharmonic Hamiltonian $H_n$.  For the sake 
of mathematical simplicity the Hamiltonians in question are given as essentially unitless operators:
\begin{eqnarray}
H_0 =& \frac{1}{2}(p^2+x^2) =& \ad a+\frac{1}{2}\\
H_n =& H_0 + \lambda x^n =& \ad a + \frac{1}{2}+ \frac{\lambda}{\sqrt{2}^n}(a+\ad)^n. \nonumber
\end{eqnarray}

As one may infer, we have made the following assumptions and simplifications:
\begin{enumerate}  
\item $\hbar  = \omega = m = 1$,\\
\item $a = \frac{x+ip}{\sqrt{2}}, \ad = \frac{x-ip}{\sqrt{2}}$, 
\item $x= \frac{a+\ad}{\sqrt{2}}, p = \frac{\ad-a}{i\sqrt{2}}$.
\end{enumerate}

Let us now give the formulation of the Lie algebras.

\begin{defn}
The Lie Algebra $\mathcal{A}_n^{(k)} = \{L_m\}_{m\in I}$ is generated by the elements
\begin{eqnarray*}
L_1 = H_0 ,& L_2=H_n 
\end{eqnarray*}

and other $L_m$ satisfying
\begin{equation}
[L_i,L_j] = \sum c_{ijm}L_m.
\end{equation}
for some structure constants $c_{ijm}\in \C.$  Furthermore, this Lie algebra should be closed under commutators 
up to order $\lambda^k$.  In other words no $L_m$ should be of the 
form $\lambda^{k+1} (a^{\dagger t}a^s - a^{\dagger s}a^t)$ for any $s,t$.  That is, formally we require 
$\lambda^{k+1} =0$ within the Lie algebra.  
\end{defn}

In the case of this paper we will consider $\mathcal{A}_n^{(1)}$ unless otherwise explicitly stated.  
In fact, [JA] only considers Lie algebras up to order one in $\lambda$ with the exceptions of $n=1$ and $n=2$ because 
these determine harmonic oscillators and their solutions are already known.  We deal with the special technique
for solving harmonic oscillators in the appendix.

Once the algebra $\mathcal{A}_n^{(1)}$ is determined we proceed in the following way.
Suppose 
\begin{equation}
[L_1,L_2] = \sum_{k=3}^{j} c_{12k}L_k
\end{equation}
where each $L_k$ is of the form
\begin{equation}
\lambda (a^{\dagger m}a^{\ell} - a^{\dagger \ell}a^m).
\end{equation}

The symmetry of these $L_k$ is important and comes back in an important way due to the normal 
ordering procedures we have adopted.  We will see this explicitly in the computations.

We then construct a unitary element of the associated Lie group by
\begin{equation}
U = \exp(\sum_{k=3}^j\alpha_k L_k).
\end{equation}

This says that the only $L_k$ allowed in our unitary are those arising directly from the commutator 
$[L_1,L_2]$. The $\alpha_k$ are real constants which we will tune as necessary.  

Once we produce such a unitary we make a transformation from $H_0$ to $H_n$ by
\begin{equation}
U^{\dagger}H_0U = H_n - \Lambda_n.
\end{equation}
In each case $\Lambda_n$ is an operator which simply controls the pertubations of eigenvalues.  
Furthermore, by the clever choice of $U$ we will have $[U,\Lambda_n]=0 + O(\lambda^2)$.
Due to the Hadamard lemma we can produce $\Lambda_n$ by computing simple commutators.

At this stage one may write the new eigenvectors as $U^{\dagger}|j\rangle$ where $|j\rangle$ 
are the eigenvectors for the harmonic Hamiltonian with eigenvalues $j+\frac{1}{2}$. 
Therefore, up to order $\lambda^2$ our equation now reads
\begin{eqnarray}
H_n U^{\dagger}|j\rangle &=& (U^{\dagger}H_0U + \Lambda_n)U^{\dagger}|j\rangle\\
&=& U^{\dagger}H_0|j\rangle + \Lambda_n U^{\dagger}|j\rangle\nonumber\\
&=& U^{\dagger}(j+\frac{1}{2})|j\rangle + U^{\dagger}\Lambda_n|j\rangle\nonumber\\
&=& (j+\frac{1}{2} + \lambda_n)U^{\dagger}|j\rangle\nonumber.
\end{eqnarray}

In essense, depending on the form of $\Lambda_n$, we will be able to read off the 
first order perturbation eigenvalues ($\lambda_n$) of 
$H_n$ with relative ease.  

A natural question arises as to when we can solve this system explicitly.  [JA] makes a passing 
statement which we will now state as a formal theorem.

\begin{thm}
If the Lie Algebra $\mathcal{A}_n$ is closed (in all orders of $\lambda$) then we can solve the 
$n^{th}$ order anharmonic oscillator in closed form.
\end{thm}

\begin{proof}
Let $\mathcal{A}_n = \{L_k \}_{k=1}^N$ be a closed Lie algebra corresponding to the Hamiltonian $H_n$.
Then consider the general Lie group element given by
\[
U = \exp(\sum_{k=1}^N \alpha_k L_k) =: \exp(L).
\]

From the Hadamard lemma we obtain
\begin{equation}
U^{\dagger}H_0U = \sum_{k=1}^{\infty}\frac{1}{k!}[L^{(k)},H_0].
\end{equation}

Since $\mathcal{A}_n$ is closed, the commutators $[L^{(k)},H_0]$ either vanish or 
give Lie algebra elements with some periodicity.
In this way we can formally sum them in a power series.  Setting our parameters to appropriate values we 
obtain
\[
U^{\dagger}H_0 U = H_n + \mathrm{perturbations}.
\]
\end{proof}

\section{Lie Algebras up to Order One in $\lambda$}

We begin the computation of the Lie algebras by giving an important commutator relation.
\begin{equation}
[\ad a, a^{\dagger k}a^{\ell} \pm a^{\dagger \ell}a^k ] = (k-\ell)(a^{\dagger k}a^{\ell}\mp a^{\dagger \ell}a^k)
\end{equation}

This equation paired with the symmetry of Weyl binomial coefficients points us to assigning 
$\adn a^m - a^{\dagger m}a^n$ as our Lie algebra elements.  Let us begin by computing $[H_0,H_n]$.
\begin{eqnarray}
[H_0,H_n] & = & [H_0 , H_0 + \frac{\lambda}{\sqrt{2}^n}(\ad+a)^n] = \nonumber \\
\lbrack H_0 , \frac{\lambda}{\sqrt{2}^n}(\ad+a)^n \rbrack & = & \frac{\lambda}{\sqrt{2}^n}[\ad a , (\ad+a)^n]  \nonumber\\
&=& \frac{\lambda}{\sqrt{2}^n}[\ad a, \sum_{k,m}\left\{\begin{array}{c}n\\m\end{array}\right\}_k 
(a^{\dagger m-k}a^{n-m-k} + a^{\dagger n-m-k}a^{m-k})]\nonumber \\ 
&=& \frac{\lambda}{\sqrt{2}^n}\sum_{k,m} \left\{\begin{array}{c}n\\m\end{array}\right\}_k 
(2m-n) (a^{\dagger m-k}a^{n-m-k} - a^{\dagger n-m-k}a^{m-k}) 
\end{eqnarray}

We will throw away the multiplicative constants in favor of rescaling them by $\alpha_k$ in our general Lie group element.  
Therefore, the first batch of elements revealed to us are those of the form $\lambda (a^{\dagger\ell} a^m - a^{\dagger m}a^{\ell})$.  
Once we realize these, we begin commuting again with $H_0$ to find more elements of the form
$\lambda (a^{\dagger\ell} a^m + a^{\dagger m}a^{\ell})$.  In order to close $\mathcal{A}_n^{(1)}$ we also need
to add a central element $I$ to our Lie algebra.  All other commutators will involve terms with $\lambda^2$
and therefore we disregard them in $\mathcal{A}_n^{(1)}$.

\begin{rem}
Notice that no elements of the form $\adn a^n$ appear anywhere.  This is because they can be written in terms of the
number operator $N=\ad a$ which commutes with $H_0$.  
\end{rem}

 Given a general Hamiltonian $H_n$, with special exceptions, by simple combinatorial formulae one infers
the number of generators for $\mathcal{A}_n^{(1)}$ by
\begin{eqnarray}
|\mathcal{A}_{2k}^{(1)}| &=& (k+1)k + 3\\
|\mathcal{A}_{2k+1}^{(1)}|& =& (k+1)k+3  \nonumber
\end{eqnarray}

\begin{ex}
Let's take a quick look at $\mathcal{A}_6^{(1)}$.  Of course, we let $L_1=H_0$ and $L_2=H_6$. 
By our computation we know the generators arising from $[L_1,L_2]$ are as follows:\\
\[
\begin{array}{lll}
&&\\
\lambda(a^{\dagger 6}-a^6), & \lambda(a^{\dagger 5}a - \ad a^5), &\lambda(\adqu a^2 - \adsq a^4),\\
\lambda(\adqu -a^4),       & \lambda(\adcu a - \ad a^3),       &\lambda (\adsq - a^2). \\
&&
\end{array}
\]

Furthermore, commuting these with $L_1$ we obtain simlarly symmetric elements with plus signs. Finally we add in $I$
to account for commuting elements.  Notice if we commute any other elements we obtain an element
in $\mathcal{A}_6^{(2)}$ which we have formally disallowed for now.
Therefore $|\mathcal{A}_6^{(1)}| = 15 = (3+1)3 + 3$ as previously stated.
\end{ex}

\section{Explicit Computations}

In this section we will derive the first order perturbation for all anharmonic oscillators with Hamiltonians
of the form
\[
H_n = \frac{1}{2}(p^2+x^2)  + \lambda x^n.
\]
There are two distinct cases for computing first order perturbations; odd and even.  We will treat the odd case first.\\

For the sake of uniformity in our calculations we will consider Hamiltonians of the form $H_{2k-1}$ and $H_{2k}$.
\subsection{Odd Powered Potentials}
From our earlier computations of $[H_0,H_{2k-1}]$ and our prescribed form of $U$ we have
\begin{equation}
U = \exp (\lambda \sum_{\ell=0}^{m-1}\sum_{m=1}^k \alpha_{m,\ell} a^{\dagger 2m-1-\ell}a^{\ell} - a^{\dagger \ell}a^{2m-1-\ell}). 
\end{equation}

By requiring $\alpha_{m,\ell}\in \R$ we obtain $U^{\dagger} = U^{-1}$ and we may now apply the Hadamard lemma.
Letting 
\[
X:= \lambda \sum_{\ell=0}^{m-1}\sum_{m=1}^k \alpha_{m,\ell} (a^{\dagger 2m-1-\ell}a^{\ell} - a^{\dagger \ell}a^{2m-1-\ell})
\]
we have:

\begin{equation}
U^{\dagger}H_0 U = H_0 + [-X,H_0] + \frac{1}{2!}[-X,[-X,H_0]] + \cdots
\end{equation}

We notice immediately that $X$ contains a multiplicative factor of $\lambda$ and since we have $\lambda^2 = 0$ we may ignore 
all terms past $[-X,H_0]$.\\

Using previous computations and elementary properties of derivations we have 

\begin{equation}
[H_0,X] =\lambda \sum_{\ell=0}^{m-1} \sum_{m=1}^{k} \alpha_{m,\ell}(2m-1-2\ell)(a^{\dagger 2m-1-\ell}a^{\ell} + a^{\dagger \ell}a^{2m-1-\ell}).
\end{equation}

If we recognize that 
\[
x^{2k-1} = \frac{(a+\ad)^{2k-1}}{\sqrt{2}^{2k-1}} =\frac{1}{\sqrt{2}^{2k-1}} 
\sum_{m=0}^{2k-1} \sum_{j=0}^{\min\{2k-1,2k-1-m\}}\left\{\begin{array}{c}2k-1\\m\end{array}\right\}_j a^{\dagger m-j}a^{2k-1-m-j},
\]

we see that in order to produce $H_{2k-1} = H_0 + \lambda x^{2k-1}$ we need to set
\begin{equation}
\alpha_{m,\ell}(2m-2\ell-1) =\frac{1}{\sqrt{2}^{2k-1}} \left\{\begin{array}{c} 2k-1\\k-m+\ell \end{array}\right\}_{k-m}.
\end{equation}

Notice what we have done.  We have transformed $H_0$ into $H_{2k-1} + O(\lambda^2)$.  Hence, there is no perturbation 
term up to first order.

\begin{ex}
Let us compute the example of $H_1 = H_0 + \lambda x$ explicitly.  We already know that this is a shifted harmonic
oscillator, where a simple change of variables reveals the energy eigenvalues are $n+\frac{1}{2} - \frac{\lambda^2}{2}$.\\
In our case the Lie algebra simplifies slightly and is fully closed as 
\[
H_0,H_1,\lambda(\ad-a),I.
\]  

Our appropriate
unitary transformation $U$ is therefore given as

\[
U = \exp(\alpha \lambda(\ad-a)).
\]

We can compute this even more explicitly than before given that $[\ad-a,\ad]=[\ad-a,a]=-1$.\\
In this case we know 
\begin{equation}
\lbrack A,B \rbrack = \beta \in \C \implies \lbrack A,e^B\rbrack = \beta e^B
\end{equation}
for any abstract operators $A,B$.

Therefore $[a,U] =\alpha\lambda U$ and $[\ad,U]= \alpha\lambda U$
\begin{eqnarray*}
U^{\dagger}H_0U &=& U^{\dagger}\ad aU +\frac{1}{2}\\
&=& U^{\dagger}\ad U(\alpha\lambda + a)+\frac{1}{2}\\
&=& U^{\dagger}U (\ad+\alpha\lambda)(a+\alpha\lambda) + \frac{1}{2}\\
&=& \ad a + \frac{1}{2} + \alpha \lambda(\ad + a) + \alpha^2\lambda^2.
\end{eqnarray*}

Setting $\alpha = \frac{1}{\sqrt{2}}$ we derive

\begin{equation}
U^{\dagger}H_0 U = H_1 + \frac{\lambda^2}{2}.
\end{equation}
This is exactly the result we previously knew.  Notice however,  there is no $\lambda$ term.  The first 
perturbation term is order $\lambda^2$.
\end{ex}

\begin{ex}
To better see the odd powered result more explicitly, let us compute the result for $H_5$.\\
Our Lie algebra up to order one is given by
\[
\begin{array}{lll}
\lambda (a^{\dagger 5} - a^5),&\lambda (\adqu a - \ad a^4), & \lambda(\adcu a^2 - \adsq a^3),\\
\lambda (\adcu - a^3),&\lambda (\adsq a- \ad a^2),& \lambda (\ad -a),\\
\lambda (a^{\dagger 5} + a^5),&\lambda (\adqu a + \ad a^4), & \lambda(\adcu a^2 + \adsq a^3),\\
\lambda (\adcu + a^3),&\lambda (\adsq a+ \ad a^2),& \lambda (\ad +a),\\
H_0,& H_5, & I.
\end{array}
\]

Therefore $U$ is
\[
U = \exp(\lambda (\alpha_{3,0} (a^{\dagger 5}-a^5) + \alpha_{3,1}(\adqu a - \ad a^4)+ \cdots + \alpha_{1,0}(\ad-a))).
\]

Our Hamiltonian transforms as
\begin{eqnarray*}
U^{\dagger}H_0 U &=& H_0 + \lambda[H_0,\alpha_{3,0} (a^{\dagger 5}-a^5) + \alpha_{3,1}(\adqu a - \ad a^4)+ \cdots + \alpha_{1,0}(\ad-a)] + 
O(\lambda^2)\\
&=&H_0+ \lambda (5 \alpha_{3,0} (a^{\dagger 5}+a^5) +3 \alpha_{3,1}(\adqu a + \ad a^4)+ \cdots + \alpha_{1,0}(\ad+a)).
\end{eqnarray*}

Setting 
\begin{eqnarray*}
\alpha_{3,0} &=& 2^{-5/2}/5, \\
\alpha_{3,1} &=& 2^{-5/2}5/3,\\
\alpha_{3,2} &=& 2^{-5/2}10,\\
\alpha_{2,0} &=& 2^{-5/2}10/3,\\
\alpha_{2,1} &=& 2^{-5/2}30,\\
\alpha_{1,0} &=& 2^{-5/2}15,
\end{eqnarray*}
we obtain

\begin{equation}
U^{\dagger}H_0 U = H_0 + \lambda x^5 + O(\lambda^2) = H_5 + O(\lambda^2).
\end{equation}
\end{ex}

\subsection{Even Powered Potentials}

It is the goal of this section to show the first order perturbation energies for oscillators corresponding to
$H_{2k}$ are
\begin{equation}
n+\frac{1}{2} + \frac{\lambda}{2^k}(\sum_{j=0}^k j! \left\{\begin{array}{c}2k\\k\end{array}\right\}_{k-j}
\left(\begin{array}{c}n\\j\end{array}\right) ).
\end{equation}

This result is obtained rather easily utilizing the technology we have developed for odd powered potentials.
We only need to realize that the operators which are not Lie algebra elements are of the form $\adn a^n$.
This comes from the required symmetry of our generators.  Hence our unitary $U$ will appear exactly as in the 
odd powered case and the Hadamard lemma yields:
\begin{equation}
U^{\dagger}H_0 U = H_0 + \lambda x^{2k} -\frac{\lambda}{2^k}\sum_{j=0}^k  \left\{\begin{array}{c}2k\\k\end{array}\right\}_{k-j} 
a^{\dagger j}a^j + O(\lambda^2).
\end{equation}

The only thing left to pretty up our example is changing $a^{\dagger k}a^k$ into an expression of number operators.

Recall $N=\ad a$ has nonnegative integer eigenvalues given by $N|n\rangle = n|n\rangle$.  In this way any expression $f(N)$
in our perturbation expansion will give eigenvalues $f(n)$ by the functional calculus.

\begin{prop}
\begin{equation}
a^{\dagger k}a^k = k! \left(\begin{array}{c}N\\k\end{array}\right).
\end{equation}
\end{prop}

\begin{proof}
We refer back to our commutation relation $[A^n,B^m]$ from $\S$1.3.  Thus we have
\begin{eqnarray*}
a^{\dagger k}a^k &=& \ad (a^{\dagger k-1}a)a^{k-1}\\
&=& \ad (a a^{\dagger k-1} - (k-1)a^{\dagger k-2})a^{k-1}\\
&=& (\ad a) a^{\dagger k-1}a^{k-1} - (k-1)a^{\dagger k-1}a^{k-1}\\
&=& (N-(k-1))a^{\dagger k-1}a^{k-1}.
\end{eqnarray*}

Repeating this we see 
\begin{equation}
a^{\dagger k}a^k = N(N-1)\cdots (N-(k-1)) = k!\left(\begin{array}{c}N\\k\end{array}\right).
\end{equation}
\end{proof}

It is merely a matter of rearranging terms to see

\begin{equation}
U^{\dagger}H_0 U = H_{2k} -\frac{\lambda}{2^k}\sum_{j=0}^k j! \left\{\begin{array}{c}2k\\k\end{array}\right\}_{k-j} 
\left(\begin{array}{c}N\\j\end{array}\right) + O(\lambda^2).
\end{equation}

\begin{ex}
Let us look briefly at the Hamiltonian 
\[
H_4 = \ad a+\frac{1}{2} + \frac{\lambda}{4}(\ad+a)^4.
\]
This is the famous quartic which has received much attention in texts and papers.  We can check our results
against those of standard perturbation theory.

Our Lie algebra $\mathcal{A}_4^{(1)}$ is given  by
\[
\begin{array}{lll}
\lambda (\adqu-a^4),& \lambda(\adcu a-\ad a^3),&\lambda(\adsq - a^2),\\
\lambda (\adqu+a^4),& \lambda(\adcu a+\ad a^3),&\lambda(\adsq + a^2),\\
H_0,& H_4, & I.
\end{array}
\]

Our unitary is given by
\[
U = \exp(\lambda (\alpha_{2,0}(\adqu-a^4)+\alpha_{2,1}(\adcu a- \ad a^3)+\alpha_{1,0}(\ad-a))).
\]

From here we must simply crank the handle for our machine and we realize
\begin{eqnarray}
U^{\dagger}H_0 U& =& H_0 + [H_0,\lambda (\alpha_{2,0}(\adqu-a^4)+\alpha_{2,1}(\adcu a- \ad a^3)+\alpha_{1,0}(\ad-a))]\nonumber\\
&=& H_0 + \lambda (4\alpha_{2,0}(\adqu+a^4)+2\alpha_{2,1}(\adcu a+ \ad a^3)+\alpha_{1,0}(\ad+a)).
\end{eqnarray}

Setting 
\begin{eqnarray*}
\alpha_{2,0} &=& 1/16 \\
\alpha_{2,1} &=& 1/2\\
\alpha_{1,0} &=& 3/4
\end{eqnarray*}

we arrive at

\begin{eqnarray}
U^{\dagger}H_0 U& =& H_4 - \frac{\lambda}{4}(6\adsq a^2 + 12\ad a + 3)\\
&=& H_4 - \frac{3\lambda}{2}(N(N+1)) - \frac{3\lambda}{4}.\nonumber
\end{eqnarray}

If we look to the ground state we see that 
\[
E_0 = \frac{1}{2} + \frac{3\lambda}{4} + O(\lambda^2)
\]

which agrees with the standard perturbation theory.
\end{ex}

In particular the perturbed ground state of the anharmonic oscillator corresponding to
$H_{2k}$ is given by
\begin{equation}
E_0 = \frac{1}{2} + \frac{\lambda}{2^k} \left\{\begin{array}{c}2k\\k\end{array}\right\}_k
+ O(\lambda^2) = \frac{1}{2} + \frac{\lambda (2k)!}{2^{2k}k!}+O(\lambda^2).
\end{equation}

\section{Extending the Method}

The second subsidiary goal of this paper is to show several extensions to this method and invite research into
even more applications of Lie algebras into physics.

\subsection{Simple One Dimensional Corollaries}

Now that we have given explicit formulae for computing perturbation eigenvalues for potentials of the form
$\lambda x^n$ we can extend by linearity (up to order one) and immediately recover eigenvalues for 
polynomial potentials.   In fact we can extend this further to convergent power series.

\begin{ex}
Let us consider 
\[
H = \ad a + \frac{1}{2} + \lambda e^x.
\]
Of course this can be rewritten 
\[
H = \ad a+\frac{1}{2} + \lambda (\sum_k \frac{x^k}{k!}).
\]

If we notice that only the even powered potentials contribute perturbations up to first order then we will also compute
perturbations for $H = H_0 + \lambda\cosh(x)$ as well.\\

Let us compute only the ground state energy.  We have
\[
E_0 = \frac{1}{2} + \lambda(\sum_k \frac{(2k)!}{2^{2k}k!(2k)!})) = \frac{1}{2} + \lambda(\sum_k \frac{4^{-k}}{k!}) = 
\frac{1}{2} + \lambda\exp(1/4).
\]

\end{ex}

Moreoever, we can add any number of perturbation parameters and solve the system accordingly.
In particular we can essentially read off first order perturbations for Hamiltonians of the form 
\[
H = H_0 + \sum_{j=1}^n \lambda_j x^{k_j}.
\]

\subsection{Simple $N$ Dimensional Corollaries}

In extending this method, it is natural to ask whether one can tackle higher dimensional systems with a similar approach.
In our case, we certainly can attack higher dimensional problems similarly, but the construction of the Lie algebra
is different.  For the simple $N$ dimensional corollaries, we will assume our oscillator potential is not coupled
(i.e. No terms of the form $\lambda x^jy^k$) appear.
For the sake of simplicity let us go through the construction of the Lie algebras for a two dimensional oscillator.\\

Consider
\[
H_{n,m} = \ad_x a_x + \frac{1}{2} + \lambda_1 x^n + \ad_y a_y +\frac{1}{2}+\lambda_2 y^m.
\]

We will take four elements as given in our Lie algebra
\begin{eqnarray*}
H_{0,0} &=& \ad_x a_x + \frac{1}{2} + \ad_y a_y + \frac{1}{2},\\
H_{0,m} &=& H_{0,0} + \lambda_2 y^m,\\
H_{n,0} &=& H_{0,0} + \lambda_1 x^n,\\
H_{n,m} &=& H_{0,0} + \lambda_1 x^n + \lambda_2 y^m.
\end{eqnarray*}

In this way we will set up our Lie algebra as two independent oscillator Lie algebras and solve our problems from before.
Consider for example
\[
H_{1,4} = H_{0,0} + \lambda_1 x + \lambda_2 y^4.
\]

Our Lie algebra $\mathcal{A}_{1,4}^{(1,1)}$ will have the following elements
\[
\begin{array}{lll}
H_{0,0},& H_{1,0},& \lambda_1(\ad_x - a_x),\\
H_{0,4}, & H_{1,4}, & I,\\
\lambda_2(\adqu_y - a^4_y),& \lambda_2 (\adcu_ya_y - \ad_y a^3_y), & \lambda_2(\adsq_y - a^2_y),\\
\lambda_2(\adqu_y + a^4_y),& \lambda_2 (\adcu_ya_y + \ad_y a^3_y), & \lambda_2(\adsq_y + a^2_y).
\end{array}
\]

Our unitary takes the form
\[
U = \exp(\alpha \lambda_1(\ad_x-a_x) + \beta_1\lambda_2(\adqu_y-a^4_y)+ \beta_2\lambda_2(\adcu_ya_y - \ad_ya^3_y)+ 
\beta_3\lambda_2(\adsq_y - a^2_y)).
\]

Now we use the Hadamard lemma again, but taking advantage of the relations
\begin{equation}
\lbrack \ad_j , a_k \rbrack = \delta_{jk}
\end{equation}

we can completely separate $x$ variables from $y$ variables and our calculation plays out exactly as before.

For the Hamiltonian $H_{1,4}$ our perturbed ground state is
\[
E_0 = \frac{1}{2} + O(\lambda_1^2) + \frac{1}{2} + \frac{3\lambda_2}{4} + O(\lambda_2^2).
\]

Now we can use all the simple one dimensional corollaries in turn as well.  

\subsection{Higher Order Perturbations}

Since perturbation theory is meant to compute more than first order terms we seek to use this Lie algebraic method
to compute higher order terms.  Certainly, one can see that using the transformations $U$ we have
set up thus far will produce higher order terms.  One can see this if we set
\[
U = \exp(\lambda L).
\]

Our transformation becomes
\begin{equation}
U^{\dagger}H_0 U = H_0 + \lambda [H_0,L] - \frac{\lambda^2}{2} [L,[H_0,L]] \dots
\end{equation}

This approach, however, changes our Hamiltonian funadmentally.  In fact, we end up not solving any problems,
but instead creating more.  A quick trial calculation with any Hamiltonian carrying term $x^3$ or higher
will reveal that we cannot cancel certain terms arising from $[L,[L,H_0]]$.  To remove this difficulty we must 
expand our Lie algebra to include terms carrying $\lambda^k$ for whichever $k$ we should choose.
It is therefore convenient to write our unitary transformation as

\begin{equation}
U = \exp(\sum_{j=1}^k \lambda^j L_{(j)}),
\end{equation}  
where $L_{(j)}$ are Lie algebra elements arising from $j^{th}$ order comutators.

\begin{ex} 
Let us return briefly to the quartic oscillator and calculate its second order perturbation.  
Since it is well studied we may verify our results easily.

Computing commutators and commutators of commutators one will arrive at the following Lie algebra
up to order 2.
\[
\begin{array}{lll}
H_0, & H_4, & I,\\
\lambda (\adqu\pm a^4),& \lambda(\adcu a \pm\ad a^3),&\lambda(\adsq \pm a^2),\\
\lambda^2 (a^{\dagger 6} \pm a^6),& \lambda^2(\adqu a^2 \pm \adsq a^4),& \\
\lambda^2(\adqu \pm a^4),& \lambda^2 (\adcu a \pm \ad a^3), & \lambda^2(\adsq \pm a^2).
\end{array}
\]

Knowing the form of our necessary first order transformation, we add four terms 
to the exponential by

\begin{eqnarray}
U &=& \exp(\lambda(\frac{1}{16}(\adqu-a^4)+\frac{1}{2}(\adcu a- \ad a^3)+ \frac{3}{4}(\adsq - a^2))\\
& & +\lambda^2\beta_1(a^{\dagger 6} - a^6)+ \lambda^2\beta_2(\adqu a^2 - \adsq a^4)\nonumber\\
& & + \lambda^2\beta_3(\adcu a - \ad a^3) + \lambda^2\beta_4(\adsq - a^2)).\nonumber\\
U&=& \exp(\lambda L_{(1)} + \lambda^2 L_{(2)}).\nonumber
\end{eqnarray}

For simplicity, let us compute only the ground state energy given by $U^{\dagger}H_0U$.
We have
\begin{eqnarray*}
U^{\dagger}H_0 U &=& \exp(-\lambda L_{(1)} - \lambda^2 L_{(2)} )H_0 \exp(\lambda L_{(1)}+ \lambda^2 L_{(2)})\\ 
&=& H_0 + \lambda [H_0,L_{(1)}] + \lambda^2[H_0,L_{(2)}] -\frac{\lambda^2}{2}[L_{(1)},[H_0,L_{(1)}]] + O(\lambda^3)
\end{eqnarray*}

From here it is a matter of computing commutators and adjusting $\beta_1,\beta_2,\beta_3,\beta_4$ to cancel higher order terms
not given as functions of number operators.\\  

When we compute the ground state energy we are concerned only with constant
terms. Therefore, looking to our commutators we have the $\lambda^2$ term
\[
\lbrack \frac{1}{16}(\adqu-a^4)+\frac{1}{2}(\adcu a- \ad a^3)+ \frac{3}{4}(\adsq - a^2),
\frac{1}{4}(\adqu+a^4)+(\adcu a+ \ad a^3)+ \frac{3}{2}(\adsq + a^2)  \rbrack.
\]

Expanding this we are left with two terms giving constants
\begin{eqnarray*}
\frac{1}{64}[\adqu-a^4,\adqu + a^4 ] & \mathrm{and}& \frac{9}{8}[\adsq-a^2,\adsq+a^2].
\end{eqnarray*}

Our constant terms turn out to be $\frac{-2(4!)}{64}$ and $\frac{-2(2!)9}{8}$ yielding $\frac{-21}{4}$.

Finally, our ground state energy up to second order will be given by

\begin{equation}
 (U^{\dagger}H_0 U)U^{\dagger}|0\rangle = (H_4 - \frac{3\lambda}{4} - \frac{-21}{4}\frac{\lambda^2}{2}) U^{\dagger}|0\rangle,
\end{equation}
yielding
\begin{equation}
E_{0} = \frac{1}{2} +\frac{3\lambda}{4} - \frac{21\lambda^2}{8} + O(\lambda^3).
\end{equation}

Indeed, this ground state energy agrees with the standard perturbation theory.\\

For the interested reader, the correct $\beta$ parameter values are
\begin{eqnarray*}
\beta_1 = \frac{1}{48}, & \beta_2 = \frac{-9}{16},\\
\beta_3 = \frac{-9}{4}, & \beta_4 = \frac{-63}{32},
\end{eqnarray*}

and the second order equation appears as
\begin{eqnarray}
U^{\dagger}H_0 U& =& H_4 -\frac{3\lambda}{2}N(N+1) -\frac{3\lambda}{4} +\nonumber\\ 
&& 51\lambda^2 \left(\begin{array}{c}N\\3\end{array} \right)+
\frac{117\lambda^2}{2}  \left(\begin{array}{c}N\\2\end{array}\right) +
36\lambda^2 N + \frac{21\lambda^2}{8} +O(\lambda^3)
\end{eqnarray}

\end{ex}

%\subsection{Coupled Oscillators}

\section{Discussion}

One important problem among many arising from this paper and which we have yet neglected to mention is the 
representation theory the Lie algebras.
For example if we are dealing with the cases $\lambda x$,$\lambda x^2$ in the one dimensional case or any 
quadratic term in higher dimensional cases, we have a closed Lie algebra.  This may not be terribly surprising as 
a closed Lie algebra offers an exact solution, and these particular potentials are simply shifted or coupled
harmonic oscillators.  However, the Lie algebras we have constructed all contain central elements.  In the cases of 
$\mathcal{A}_1$ and $\mathcal{A}_2$ we have 4 dimensional closed Lie algebras with center.  Representation theory 
tells us that these are isomorphic to $\mathfrak{gl}_2$.  It remains to be seen exactly the relationship
between these Lie algebras and the symmetries they describe.  In this paper we have only used them
to calculate perturbed energy levels.  It is entirely possible there is a simpler approach to this problem from
an entirely representation theoretic standpoint.\\

As it stands, we have given explicit constructions for Lie algebras up to any order and the 
method by which we may construct a unitary operator to make the transformation
\[
H_0 \mapsto H_n + \mathrm{perturbations\ up\ to\ } O(\lambda^k). 
\]

By taking advantage of our symmetric construction of these Lie algebras, the Hadamard lemma, and several
formulae concerning abstract Weyl algebras we have managed to give eigenvalues in agreement with standard methods.\\

Another issue which we have neglected to resolve it how to deal with coupled oscillators in general.  
In the appendix we briefly mention the way to deal with potentials of the form $\lambda xy$.  The Lie algebra 
computations for coupling terms of the form $\lambda x^ny^m$ are more taxing and trickier.  This method, however, 
should be able to deal with situation, but some coordinate change may be required first.

It is the hope of the author that anharmonic oscillators are simply a useful class of examples for the propitiation of this method.
Furthermore, it is hoped that this method will help to give rise to additional representation theoretic methods in 
physics.  

\section*{Appendix: Dealing with Harmonic Oscillators}
Two important cases we haven't touched upon are those of actual harmonic oscillators in one and multiple dimensions.
Consider for example the two Hamiltonians
\begin{eqnarray}
H_2& =& \ad a + \frac{1}{2} + \lambda(\frac{\ad+a}{\sqrt{2}})^2\\
H_c &=& \ad_x a_x + \ad_ya_y +1 + \frac{\lambda}{2}(\ad_x+a_x)(\ad_y+a_y).\nonumber
\end{eqnarray}

These correspond to the shifted frequency oscillator in one dimension with new frequency $\sqrt{1+2\lambda}$
and a coupled oscillator in two dimensions with quadratic coupling term.  These cases have been well studied
and so we have neglected them thus far.  However, the techniques to compute the perturbations are special
because these Hamiltonians along with $H_0$ and $H_{0,0}$ produce closed Lie algebras.  By our theorem earlier
we know that we can solve these exactly and not concern ourselves with $k^{th}$ order perturbations.\\

The main technique we employ is to transform our ladder operators via the so called Bogoliubov tranforms.
In one dimension we have
\begin{eqnarray}
b^{\dagger}&=& U^{\dagger}\ad U = \sigma\ad + \tau a\\
b &=& U^{\dagger}a U = \sigma a+ \tau\ad.\nonumber
\end{eqnarray}

In this way we produce the new Hamiltonian
\[
\sqrt{1+2\lambda}(b^{\dagger}b + \frac{1}{2}) = \ad a + \frac{1}{2} + \frac{\lambda}{2}(\ad+a)^2.
\]
This algebra to move from $\sigma,\tau$ to this clean form of the new Hamiltonian is tedious to be sure.
The interested reader should confer with [JA] or email the author for a small set of notes.\\

A similar technique can be used for the quadratic coupling, but the transformation must take into account 
much more coupling.  Our transformation should look something like
\begin{equation}
\left(\begin{array}{c}b^{\dagger}_x\\b_x\\b^{\dagger}_y\\b_y\end{array}\right) = 
\left(\begin{array}{cccc}
\alpha & \beta & \gamma & \delta\\
\beta & \alpha & \delta & \gamma\\
\sigma & \tau & \mu & \nu\\
\tau & \sigma & \nu & \mu
\end{array}\right)
\left(\begin{array}{c}a^{\dagger}_x\\a_x\\a^{\dagger}_y\\a_y\end{array}\right)
\end{equation}

This is simply a coordinate change which decouples the coordinate variables.  The matrix, however, will take a very special form
so that $[b_i,b^{\dagger}_j] = \delta_{ij}$ as did the initial coordinates.  

\begin{rem}
Notice here that we can couple our ladder operators in many more ways. 
For example we can tackle problems such as dynamic coupling 
\[
H = x^2-\partial^2_x + y^2 - \partial_y^2 + \lambda\partial_x\partial_y,
\]

or oscillators in a magnetic field 
\[
H = x^2-\partial^2_x + y^2 - \partial_y^2 + \lambda(y\partial_x - x\partial_y).
\]
\end{rem}

So long as our coupling term contains terms of order two or less in each of the ladder operators we can 
tackle thess problems with a simple coordinate change.

\end{document}